\documentstyle[amssymb]{amsart}
\newcommand{\stopproof}{\hfill \nobreak\medskip $\blacksquare$ \\
\hspace*{\fill}}
\newcommand{\dom}{\mbox{\rm dom}}
\newcommand{\cof}{\mbox{\rm cof}}

\newcommand{\AND}{\mbox{ \rm and }}

\newcommand{\IF}{\mbox{ \rm if }}

\newcommand{\forces}[2]{\Vdash_{#1} \mbox{``} #2 \mbox{''}}

\newcommand{\proof}{{\bf Proof:} \ }

\newcommand{\PP}{{\Bbb P}}

\newcommand{\CC}{{\Bbb C}}

\newcommand{\QQ}{{\Bbb Q}}

\newcommand{\presup}[2]{\, ^{#1} \! #2}
\newcommand{\fomom}{\presup{\omega}{\omega}}

\newtheorem{lemma}{Lemma}[section]

\title{Addendum to ``Maximal Chains in $\fomom$ and Ultrapowers
of the Integers''  }
\author{Saharon Shelah}
\address{Institute of Mathematics \\ Hebrew University \\ Jerusalem,
Givat Ram,  Israel and Department of Mathematics\\
Rutgers University\\ New Brunswick, New Jersey}
\author{Juris Stepr\={a}ns}
\address{Department of Mathematics, York University \\
4700 Keele Street \\ North York, Ontario \\ Canada \ \ \ \ M3J 1P3}
\begin{document}
\maketitle

This note is intended as a supplement and clarification
to the proof of Theorem 3.3 of
 \cite{original}; namely, it is consistent that ${\frak b} = \aleph_1$
yet for every ultrafilter $U$ on $\omega$ there
 is a $\leq^*$ chain
$\{f_\xi : \xi\in\omega_2\}$ such that $\{f_\xi/U :\xi\in\omega_2\}$ is
cofinal in $\omega/U$.

 The general outline of the the proof remains the same. In other
words, a ground model is taken which satisfies
$2^{\aleph_0}=\aleph_1$ and in which there is a $\lozenge_{\omega_2}$ sequence
$\{D_\xi :\xi\in\omega_2\}$ such that for every $X\subseteq \omega_2$
there is a stationary set of ordinals, $\mu$, such that $\cof (\mu )=
\omega_1$ and such that $X\cap \mu = D_\mu$. Actually, a coding will
be used to associate with subsets of $\omega_2$, names  for subsets
on $\omega$ in certain partial orders.
 The details of this coding will be ignored except to state that
$c(D_\eta)$ will denote
the coded set and that if $\PP_{\omega_2}=\lim\{\PP_\xi : \xi\in\omega_2\}$
 is the finite support iteration of ccc
partial orders of size no greater than $\omega_1$ and
$1\forces{\PP_{\omega_2}}{X\subseteq [\omega]^\aleph_0}$ then there is a
stationary set  $S_X\subseteq\omega_2$, consisting of ordinals of
uncountable cofinality, such that
$1\forces{\PP_\xi}{X\restriction\PP_\xi = c(D_\xi)}$ for each $\xi
\in S_X$. Here $X\restriction\PP_\xi$ denotes the $\PP_\xi$-name
obtained by considering only those parts of $X$ that mention
conditions in $\PP_\xi$; to be more precise here would requires
providing the details of a specific development of names in the
theory of forcing, and so this will not be done.

The partial order $\PP_{\omega_2}$ is defined by induction using the
$\lozenge_{\omega_2}$ sequence. Simultaneously,
a partial ordering $\prec$ will be defined on
$\omega_2$ by $\eta \prec \zeta$ if and only if
\begin{itemize}
\item  $1\forces{\PP_\eta}{c(D_\eta) \mbox{ is an ultrafilter on } \omega}$
\item  $1\forces{\PP_\zeta}{c(D_\zeta) \mbox{ is an ultrafilter on } \omega}$
\item $1\forces{\PP\zeta}{c(D_\eta) = c(D_\zeta)\cap V^{\PP_\eta}}$
\end{itemize}
If $\alpha \in \omega_2$ then the order type of $\{\beta \in \alpha :
\beta \prec \alpha\}$ will be denoted by $o(\alpha)$.
Furthermore, an enumeration
 $\{g_\xi : \xi\in\omega_2\}$ will be constructed by induction
along with $\PP_{\omega_2}$
which will list all $\PP_{\omega_2}$-names for functions from
$\omega$ to $\omega$.

If $\PP_\xi$ has been defined and $1\forces{\PP_\xi}{c(D_\xi)\mbox{
is an ultrafilter on }\omega}$ then
$\PP_{\xi + 1}$ is defined to be $\PP_\xi * \CC_\xi * \QQ_\xi$ where
$\CC_\xi$ is simply Cohen forcing which adds a single generic function
$A_\xi : \omega \to 2$ and $\QQ_\xi$ adds a function $F_\xi:\omega
\to\omega$ such that $F_\xi \geq^* F_\mu$ and $F_\xi\restriction
A_\mu^{-1}\{k\} \geq^* g_{o(\mu)}\restriction A_\mu^{-1}\{k
\}$, for a certain $k\in 2$,
 for all $\mu$ such that
$\mu\prec \xi$. To be more precise, $\QQ_\xi$ is defined, in the
forcing extension by $\PP_\xi$, to consist of all pairs
$(f,\Delta)$ such that $f$ is a finite partial function from $\omega$
to $\omega$ and $\Delta \in [\xi]^{<\aleph_0}$, and the ordering is
defined by $(f,\Delta) \leq (f',\Delta')$ if
\begin{itemize}
    \item
 $f\subseteq f'$ \item
$\Delta\subseteq \Delta'$ \item if $\mu \in \Delta$ and $\mu \prec \xi$
and $m\in \dom(f'\setminus f)$ then $f'(m) \geq F_\mu(m)$
\item if $\mu \in \Delta$, $\mu \prec \xi$,
$m\in \dom(f'\setminus f)$, $1\forces{\PP_\xi}{A_\mu^{-1}\{k\}\in c(D_\xi)}$
 and $A_\mu^{-1}(m) = k$ then $f'(m) \geq g_{o(\mu)}(m)$
\end{itemize}
If  $1\forces{\PP_\xi}{c(D_\xi)\mbox{
is an ultrafilter on }\omega}$ fails to be true then
 $\QQ_\xi$ is defined to be empty. At limits the iteration is with
finite support.

To see that for every ultrafilter on $\omega$ there is an
increasing $\leq^*$ chain which is cofinal in the ultrapower, let
$G$ be $\PP_{\omega_2}$ generic over $V$ and let $\cal U$ be an
ultrafilter on $\omega$ in $V[G]$. There must be some name $U$ such that
$1\forces{\PP_{\omega_2}}{U\mbox { is an ultrafilter on }\omega}$ and
$\cal U$ is the interpretation $U$ in $V[G]$. It is well known that
there is a set which is closed under increasing $\omega_1$ sequences,
$C$ such that $1\forces{\PP_\xi}{U\mbox { is an ultrafilter on
}\omega}$  for each $\xi\in C$. It follows that if $\alpha \in \beta$ and
$\{\alpha,\beta\}\subseteq C\cap S_U$ then $\alpha\prec \beta$. It
is now easy to verify that $\{F_\xi : \xi\in C\cap S_U\}$ is an
$\leq^*$-increasing sequence. Moreover, because $C\cap S_U$ is
cofinal in $\omega_2$ it follows that $\{o(\xi) : \xi \in C\cap S_U\}
= \omega_2$ and hence $1\forces{\PP_{\omega_2}}
{\{g_{o(\xi)} : \xi \in C\cap S_U\} = \fomom}$. Therefore, if
$1\forces{\PP_{\omega_2}}{g:\omega\to\omega}$ there is some $\xi\in
C\cap S_U$ such
that $1\forces{\PP_{o(\xi)}}{(g = g_{o(\xi)} }$ and so it follows that
$$1\forces{\PP_{\omega_2}}{(\forall^\infty n\in A_\xi^{-1}\{k\})
F_\eta(n)\geq g_{o(\xi)}(n)\AND A_\xi^{-1}\{k\}\in
c(D_\eta)\subseteq U}$$ for any  $\eta\in C\cap S_U\setminus \xi$ It
follows immediately that $\{F_\xi : \xi \in C\cap S_U\}$ is cofinal
in the ultrapower by $\cal U$.

The only thing which now has to be proved  is that $\PP_{\omega_2}$ is
locally Cohen since this immediately implies that ${\frak b} = \aleph_1$.
A condition $p\in\PP_{\omega_2}$ will be said to be determined if
there is some $\Sigma_p \in [\omega_2]^{<\aleph_0}$ such that
 $\Sigma_p $ is the support of $p$
and for each $\sigma \in \Sigma_p$ there is a quadruple
$(a^\sigma_p,f^\sigma_p, \Delta^\sigma_p,g^\sigma_p)$ such that:
\begin{itemize}
\item $p\restriction \sigma\forces{\PP_\sigma}{p(\sigma)=
a^\sigma_p*(f^\sigma_p, \Delta^\sigma_p)}$ for each $\sigma \in \Sigma_p$
\item $\Delta^\sigma_p\subseteq \Sigma_p\cap \sigma$ for each $\sigma
\in \Sigma_p$
\item $p\restriction
\sigma\forces{\PP_\sigma}{g_{o(\sigma)}\restriction\dom(a^\sigma_p) =
g^\sigma_p}$ for each $\sigma\in \Sigma_p$
\item for each $\{\sigma,\tau\}\in [\Sigma_p]^2$ such that $\sigma
\prec \tau$ there is some $k_p(\sigma,\tau)\in 2$ such that
$p\restriction\tau
\forces{\PP_\tau}{A_\sigma^{-1}\{k_p(\sigma,\tau)\}\in D_\tau}$
\item $\dom(f^\sigma_p)\supseteq \dom(a^\sigma_p)$ for each
$\sigma\in \Sigma_p$
\item $\dom(f^\tau_p) \subseteq \dom(f^\sigma_p)$ for each
$\{\sigma,\tau\}\in [\Sigma_p]^2$ such that $\sigma\prec \tau$
\end{itemize}
This definition of determined differs in a substantial way
 from the definition of
{\em somewhat determined} in \cite{original}. The next lemma shows that
every condition can be extended to a determined condition; this is
problematic for the somewhat determined conditions.
\begin{lemma}
    The set of\label{determined} determined conditions is dense in
$\PP_{\omega_2}$.
\end{lemma}
\proof Induction on $\alpha\in \omega_2+1$
will be used to prove the following stronger statment:
For each $m\in\omega$ and each $p\in\PP_{\alpha}$ there is a
determined condition $q\geq p$ such that if $\sigma$ is the maximal
element of $\Sigma_q$ then $m\subseteq a^\sigma_q$ and $\sigma$ is
the maximal element of the support of $p$. Note that $a^\sigma_q$ has
the smallest domain of any function appearing in $q$ so the requirement
that $m\subseteq a^\sigma_q$ implies that $m$ is in the domain of any function
appearing in $q$.

To prove this, suppose the statement is true for all $\alpha \in
\beta$. If $\beta $ is a limit ordinal the result follows from the
finite support of the iteration; therefore suppose that $\beta =
\gamma + 1$. Then extend $p$ so that
$p\forces{\PP_{\gamma}}{p(\gamma) = a*(f,\Delta)}$. By extending, it
may be assumed that $m\subseteq
\dom (a)\subseteq \dom(f)$. Let $\bar{m}$ be the maximal element of
$\dom(f)$. Let $p'\geq p\restriction \gamma$ be such that $\Delta$ is
contained in the support of $p'$.

There are now two cases to consider: Either $\beta$ is a successor in
$\prec$ or it is a limit. If it is a successor then let $\beta^*$ be
the predecessor of $\beta$ in $\prec$. Otherwise, let $\beta^*$ be
such that $\beta^*$ is greater then the support of of $p'$ and
$\beta^* \prec \beta$ and $\beta^*$ is the successor of $\beta^{**}$ in
the ordering $\prec$. In the first case, let $p''\geq p'$ be such
that $p''\forces{\PP_\gamma}{A_{\beta^{*}}^{-1}{k}\in D_{\beta}}$. In
the second case,
choose $p''$ such that
 $p''\forces{\PP_{\beta^*}}{A_{\beta^{**}}^{-1}{k}\in D_{\beta^*}}$
and such that  $\beta^{**}$ belongs to the support of $p''$.

Now use the induction hypothesis to find a determined condition $q$
such that if $\sigma$ is the maximal element of $\Sigma_q$ then $\bar{m}\in
\dom(a^\sigma_q)$. Moreover, in the case that $\beta$ is a limit of
$\prec$, then the induction hypothesis can be used to ensure that
 $\sigma < \beta^*$. It will be shown  that the transitivity of $ \prec$
guarantees that $q*p(\gamma) = r$ is a determined condition satisfying
the extra induction requirements. Let $\Sigma_r = \Sigma_q\cup\{\beta\}$
and let $f^\sigma_r$, $a^\sigma_r$ and $\Delta^\sigma_r$ have the
values inherited from $q$ and $p(\beta)$. Furthermore,
$k_{r}(\alpha,\tau)$ can be defined to be $k_q(\alpha,\tau)$ unless
$\beta=\tau$. Here the choice of $p''$ helps.

In the case that $\beta$ is the successor of $\beta^*$, then $p''$
decides that $A_{\beta^{*}}^{-1}{k}\in D_{\beta^*}$ so
$k_{r}(\beta^*,\beta)$ can be defined to be $k$ and, moreover
$k_{r}(\mu,\beta)$ can be defined to be $k$ for each $\mu \in
\Sigma_q$ such that $\mu\prec \beta^*$. Since $\beta$ is the successor
of $\beta^*$ in $\prec$ there are no new instances with which to deal.
In the case that $\beta$ is a limit in the partial order $\prec$,
it is possibe to define $k_{r}(\beta^{**},\beta) = k$ because of the
transitivity of $\prec$. For the same reason it is possible to define
$k_{r}(\mu,\beta)$ to be $k$ for each $\mu \in \Sigma_q$ such that
$\mu\prec \beta^{**}$.  Since the support of $q$ is contained in
$\beta^*$ and $\beta^*$ is the successor of $\beta^{**}$ in the
partial order $\prec$, it follows that there are no new instances to
consider in this case as well.  \stopproof

\begin{lemma}
    The partial order $\PP_{\omega_2}$ is locally Cohen.
\end{lemma}
\proof
Let $X\in [\PP_{\omega_2}]^{\aleph_0}$. Let ${\frak M}$ be a
countable elementary submodel of $H(\omega_3)$ which contains $X$ and
the $\lozenge$-sequence $\{D_\xi : \xi\in\omega_2\}$ as well as
$\PP_{\omega_2}$.  It suffices to
show that if $p\in\PP_{\omega_2}$ and $D\subseteq {\frak M}\cap
\PP_{\omega_2}$ is a dense subset of the partial order ${\frak M}\cap
\PP_{\omega_2}$ then there is $q\in D$ and $r\in\PP_{\omega_2}$ such
that $r \geq p$ and $r\geq q$.

Given $p\in \PP_{\omega_2}$, by using Lemma~\ref{determined}, it may,
without loss of generality, be assumed that $p$ is determined.
Using the elementarity
of ${\frak M}$ it follows that there is some determined condition
$p'$ which is isomorphic to $p$. In particular, there is an order
preserving bijection $I:\Sigma_p\to \Sigma_{p'}$ such that $I$ is
the identity on $\Sigma{p}\cap \frak M$, $a_p^{\sigma} = a_{p'}^{I(\sigma)}$,
 $f_p^{\sigma} = f_{p'}^{I(\sigma)}$,
 $g^{\sigma}_p = g^{I(\sigma)}$
 and  $I$ preserve the partial ordering $\prec$.
It is not required that $\Delta^\sigma_p = \Delta^{I(\sigma)}_{p'}$
because $\Delta_{p'}^\sigma$ will be defined to be $\Sigma_{p'}\cap\sigma$.

Now let $q\in D$ be a condition extending $p'$. Using
Lemma~\ref{determined} it may again be assumed that $q$ is determined.
 It must be shown how
to define $r\in
\PP_{\omega_2}$ extending both $q$ and $p$. In order to do this,
define $s(\alpha)$ to be the unique, minimal ordinal  $\delta \in{\frak M}$
such that $\alpha\prec\delta$ if such a unique ordinal exists. Notice
that if $\alpha\notin \frak M$ and there is some $\delta\in \frak M$
such that $\alpha\prec\delta$ then $s(\alpha)$ exists. The reason for
this is that the only way that $s(\alpha)$ can fail to exist in this
context is that there are two minimal ordinals $\delta\in \frak M$
 and $\delta'\in \frak M$
such that $\alpha\prec \delta$ and $\alpha\prec \delta'$. However,
this means that the supremum of $\{\gamma : \gamma\prec \delta \AND
\gamma\prec \delta'\}$ belongs to $\frak M$ and hence there is some
$\alpha'\in {\frak M}\setminus \alpha$ such that  $\alpha'\prec \delta $ and
$\alpha'\prec \delta'$. From the easily verified fact that $\prec$ is a
tree ordering it follows
$\alpha\prec \alpha'$ contradicting the minimality assumption on
$\delta$ and $\delta'$.

Now define $r$ as follows:
\begin{itemize}
\item the domain of $r$ is the union of the domains of $q$ and $p$
\item if $\alpha \in {\frak M}$ then $r(\alpha) = a^\alpha_q * (f^\alpha_q,
\Delta^\alpha_q\cup \Delta^\alpha_p)$
\item if $\alpha\notin\frak M$ and there does not exist $\delta\in
\dom(q)$ such that $\alpha \prec \delta$ then $r(\alpha) = p(\alpha)$
\item  if $\alpha\notin\frak M$ and there exists $\delta\in
\dom(q)$ such that $\alpha \prec \delta$ then recall that $s(\alpha)$ is
defined and define $r(\alpha) =
a^\alpha_r*(f^\alpha_r,\Delta^\alpha_p)$  where the function
$a^\alpha_r$ is defined by
$$ a^\alpha_r(n) = \left\{\begin{array}{ll}
                          a_p^\alpha(n)    &  \IF n\in\dom(a_p^\alpha)  \\
                          k_p(\alpha,s(\alpha)) + 1\mod 2    &  \IF
 n\notin\dom(a_p^\alpha)
                          \end{array}\right.$$
(note that in this case $k_p(\alpha,s(\alpha))$ has a natural definition
because  $\prec$ is a tree ordering) and the function $f^\alpha_r$
is defined by
$$ f^\alpha_r(n) = \left\{\begin{array}{ll}
                          f_p^\alpha(n)    &  \IF n\in\dom(f_p^\alpha)  \\
\min\{f_q^\beta(n) : \beta\in \Sigma_p\cap{\frak M}\AND
\alpha\prec\beta\}  &  \IF  n\notin\dom(a_p^\alpha)
                          \end{array}\right.$$
\end{itemize}

The fact that $r\geq q$ is immediate because $a_r^\mu =
a_q^\mu$ and $f_r^\mu = f_q^\mu$ for each $\mu\in \Sigma_q$
and, moreover, if $\alpha \in\dom(q)$ then $\Delta_q^\alpha\subseteq \frak
M$; so there is no restriction on the points in the domain of $r$ not
in the domain of $q$.

It will be shown that $r \geq p$ by inductively proving that $r\restriction
\rho \geq p\restriction \rho$ for each $\rho \in\omega_2$. If $\rho =
0$ there is nothing to do and at limits the finite support of the
iteration makes the task easy. So suppose that $r\restriction
\rho \geq p\restriction \rho$.
Tt suffices to  show that the following {\em Key Condition}
is satisfied: If
\begin{itemize}
\item $\alpha \prec \beta\leq \rho$
\item $\beta\in \Sigma_p$
\item $\alpha \in \Delta_p^\beta$
\item $n$ is in the domain of $f_r^\beta\setminus f_p^\beta$
\end{itemize} then
$f_r^\beta(n)\geq f_r^\alpha(n)$ and, in addition, if $a_r^\alpha(n) =
k_r(\alpha,\beta)$  then
$r\restriction (\rho + 1)
\forces{\PP_{\rho + 1}}{f^\beta_r(n)\geq g_{o(\alpha)}(n)}$. This
will be established by considering various cases.

{\noindent \bf Case 1}\\
Suppose that $\alpha$ and $\beta$ both belong to $\frak M$. Since
$q\geq p'$, from the definition of $p'$ and the partial order $\QQ_\beta$
 it easily follows
that the Key Condition
is satisfied. There is no need to use the induction hypothesis in
this case.

{\noindent \bf Case 2}\\
Suppose now that $\beta$ belongs to $\frak M$ but $\alpha$ does not.
First it will be shown that $f^\beta_r(n)\geq f^\alpha_r(n)$.
There are two subcases to consider; either $n$ belongs to the domain
of $f_p^\alpha$ or it does not. If it does, then
$f^{I(\alpha)}_{p'}(n) = f^\alpha_p(n) = f^\alpha_r(n)$ and $I(\alpha) \in
\Delta_{p'}^\beta$. Because $I(\alpha) \prec \beta$, it follows from
the fact that $q\geq p'$ that $f^\beta_q(n) \geq
f^{I(\alpha)}_{p'}(n) = f^\alpha_r(n)$. The other possibility is
that $n$ does not belong to the domain of $f_p^\alpha$. In this case,
the definition of $r$ asserts that
$$f^\alpha_r(n) = \min\{f_q^\gamma(n) : \alpha \prec \gamma \AND
\gamma \in {\frak M}\cap\dom(p)\}$$
and, since $\beta$ is in the support of $p$ and $\alpha \prec \beta$,
it follows that $f^\alpha_r(n) \leq f^\beta_q(n)$.

It must now be shown that, if $k_r(\alpha,\beta) = a_r(n)$ then
$r\restriction (\rho + 1)\forces{\PP_{\rho + 1}}
{f^\beta_r(n) \geq g_{o(\alpha)}(n)}$.
There are again two subcases to consider; either $n$ belongs to the domain
of $a_p^\alpha$ or it does not. If it does, then
$g^{I(\alpha)}_{p'}(n) = g^\alpha_p(n)$ and $I(\alpha) \in
\Delta_{p'}^\beta$. Because $I(\alpha) \prec \beta$, it follows from
the fact that $q\geq p'$ that $q\restriction\beta\forces{\PP_\beta}
{f^\beta_q(n) \geq
g^{I(\alpha)}_{p'}(n)}$ while, on the other hand,
$p\restriction\alpha\forces{\PP_\alpha}{g_{o(\alpha)}(n) = g^\alpha_p(n)}$.
Since the induction hypothesis implies that $r\restriction\alpha \geq
p\restriction\alpha$ and it has already been noted that $r\geq q$ it
follows that
$r\forces{\rho + 1}{f^\beta_q(n) \geq
g^{I(\alpha)}_{p'}(n) =   g^\alpha_p(n) = g_{o(\alpha)}(n)}$.
 The other possibility is
that $n$ does not belong to the domain of $a_p^\alpha$. In this case,
the definition of $r$ guarantees that $a_r^\alpha(n)\neq
k_r(\alpha,s(\alpha))$ and, the minimality of $s(\alpha)$ guarantees
that $s(\alpha) \prec \beta$ because $\alpha\prec \beta$ and $\beta
\in \frak M$. The transitivity of $\prec$, now guarantees that
$k_r(\alpha,\beta) \neq a_r^\alpha(n)$ and so the Key Condition is
vacuously satsified.

{\noindent \bf Case 3}\\
Suppose now that $\alpha\in \frak M$ but $\beta\notin \frak M$.
It will first be shown that $f^\beta_r(n) \geq f^\alpha_r(n)$. To see
this, recall that $f^\beta_r(n) = f^\gamma_q(n)$ for some $\gamma$
such that $\beta\prec\gamma$, $\gamma\in \frak M$ and $\gamma$
belongs to the support of $p$ --- recall that it is being assumed
that $n$ is not in the domain of $f^\beta_p$ and this function was
only extended in the case that there was an appropriate $\gamma$.
Recall also that this implies that
$\Delta_{p'}^\gamma = \Sigma_{p'}\cap \gamma\cap\frak M$ and hence
$\alpha\in \Delta_{p'}^\gamma$. Because
 $\alpha \prec \beta\prec\gamma$ it follows from the fact that $q\geq
p'$ that $f^\alpha_q(n) \leq f^\gamma_q(n) = f^\beta_r(n)$.

Now consider $g_{o(\alpha)}(n)$. Since $\alpha \in
\Delta_{p'}^\gamma$ and $\alpha\prec \gamma$ it follows that
$$q\forces{}{g_{o(\alpha)}\leq f_q^\gamma(n)}$$ and, because it has
already been noted that $r\geq q$ it follows that
$$r\restriction (\rho + 1)\forces{\PP_{\rho + 1}}{g_{o(\alpha)}\leq f_q^\gamma(n)}$$
Since $f^\gamma_q(n) = f^\beta_r(n)$ it follows that the Key
Condition has been satisfied.

{\noindent \bf Case 4}\\
Finally, suppose that neither $\alpha$ nor $\beta$ belongs to $\frak M$.
To show that $f_r^\beta(n) \geq f_r^\alpha(n)$ two cases must again
be considered; either $n$ belongs to the domain of $f_p^\alpha$ or it
does not. If it does, then $f^\beta_r(n) = f^\gamma_q(n)$ for some $\gamma $
such that $\beta\prec\gamma$, $\gamma\in \frak M$ and $\gamma$
belongs to the support of $p$. Since $\alpha \prec \beta \prec
\gamma$ it follows that $I(\alpha) \prec \gamma$ and so
$f_q^\gamma(n) \geq f_{p'}^{I(\alpha)}(n) = f_p^\alpha(n) = f_r^\alpha(n)$.
On the other hand, if $n$ does not belong to the domain of $f_p^\alpha$ then
 $$f^\alpha_r(n) = \min\{f^\gamma_q(n) :  \alpha\prec\gamma \AND
\gamma\in \frak M\cap \dom(p)\}$$ and, since this minimum is taken
over a set which includes $\gamma$, it follows that $f^\alpha_r(n)\leq
f^\gamma_q(n) = f^\beta_r(n)$.

To show that $r\restriction (\rho + 1)\forces{\PP_{\rho + 1}}
{g_{o(\alpha)}(n) \leq f^\beta_r(n)}$ there
are, once again, two cases to consider; either $n$ belongs to the domain
of $a_p^\alpha$ or it does not. If it does, then
$g^{I(\alpha)}_{p'}(n) = g^\alpha_r(n)$ and $I(\alpha) \in
\Delta_{p'}^\gamma$. Because $I(\alpha) \prec \gamma$, it follows from
the fact that $q\geq p'$ that $q\forces{}{f^\gamma_q(n)\geq
g^{I(\alpha)}_{p'}(n) }$. On the other hand,
$p\restriction\alpha \forces{\PP_\alpha}{g_{o(\alpha)}(n) =
g^\alpha_r(n)}$ and so the, because the induction hypothesis yields that
$r\restriction \alpha \geq p\restriction \alpha$ it follows that
$r\restriction(\rho  + 1) \forces{\PP_{\rho + 1}}
{f^\beta_q(n) = f^\gamma_q(n)\geq  g^{I(\alpha)}_{p'}(n)
= g^\alpha_r(n) = g{o(\alpha)}(n)}$.
 The other possibility is
that $n$ does not belong to the domain of $a_p^\alpha$. In this case,
the definition of $r$ guarantees that $a_r^\alpha(n)\neq
k_r(\alpha,s(\alpha))$. The fact that $\alpha \prec \beta$, together
with the uniquenness of $s(\alpha)$ guarantees
that $s(\alpha) = s(\beta)$.  The transitivity of $\prec$, now guarantees that
$k_r(\alpha,\beta) = k_r(\alpha,s(\beta)) = k_r(\alpha,s(\alpha))
 \neq a_r^\alpha(n)$ and so the Key Condition is
vacuously satisified. The use of $ k_r(\alpha,s(\beta))$ and
$ k_r(\alpha,s(\alpha))$ here is a slight abuse of notation because
there is no guarantee that $s(\alpha)$ belongs to the domain of $r$.
Nevertherless, because $k_r(\alpha,\gamma) $ is defined for some
$\gamma$ such that $\alpha \prec s(\alpha) \prec \gamma$ there is no
harm in this abuse.
\stopproof
\makeatletter \renewcommand{\@biblabel}[1]{\hfill#1.}\makeatother
\renewcommand{\bysame}{\leavevmode\hbox to3em{\hrulefill}\,}


\begin{thebibliography}{1}

\bibitem{original}
S.~Shelah and J.~Stepr\={a}ns, {\em Maximal Chains in $\fomom$ and Ultrapowers
of the Integers}, to appear in Arch. f\"{u}r Math. Log.

\end{thebibliography}
\end{document}